\documentclass[]{interact}

\usepackage{epstopdf}
\usepackage[caption=false]{subfig}

\usepackage[numbers,sort&compress]{natbib}
\bibpunct[, ]{[}{]}{,}{n}{,}{,}
\usepackage{algorithmicx,algorithm}

\usepackage{amsmath}

\usepackage[noend]{algpseudocode}

\usepackage{algorithmicx,algorithm}

\usepackage[colorlinks, linkcolor=blue, anchorcolor=blue, citecolor=blue]{hyperref}

\theoremstyle{plain}
\newtheorem{theorem}{Theorem}
\newtheorem{lemma}{Lemma}
\newtheorem{corollary}{Corollary}
\newtheorem{proposition}{Proposition}

\theoremstyle{definition}
\newtheorem{definition}{Definition}

\theoremstyle{remark}
\newtheorem{remark}{Remark}

\begin{document}


\title{High-Dimensional Dynamic Systems Identification with Additional Constraints}

\author{
\name{Junlin Li \textsuperscript{\dag}\thanks{\textsuperscript{\dag} Huazhong University of Science and Technology, Wuhan 430074, Hubei, P. R. China}
}
}


\maketitle

\begin{abstract}
This note presents a unified analysis of the identification of dynamical systems with low-rank constraints under high-dimensional scaling. This identification problem for dynamic systems are challenging due to the intrinsic dependency of the data. To alleviate this problem, we first formulate this identification problem into a multivariate linear regression problem with row-sub-Gaussian measurement matrix using the more general input designs and the independent repeated sampling schemes. We then propose a nuclear norm heuristic method that estimates the parameter matrix of dynamic system from a few input-state data samples. Based on this, we can extend the existing results. In this paper, we consider two scenarios. \emph{(i)} In the noiseless scenario, nuclear-norm minimization is introduced for promoting low-rank. We define the notion of \emph{weak restricted isometry property}, which is weaker than the ordinary restricted isometry property, and show it holds with high probability for the row-sub-Gaussian measurement matrix. Thereby, the rank-minimization matrix can be exactly recovered from finite number of data samples. \emph{(ii)} In the noisy scenario, a regularized framework involving nuclear norm penalty is established. We give the notion of \emph{operator norm curvature condition} for the loss function, and show it holds for row-sub-Gaussian measurement matrix with high probability. Consequently, when specifying the suitable choice of the regularization parameter, the operator norm error of the optimal solution of this program has a sharp bound given a finite amount of data samples. This operator norm error bound is stronger than the ordinary Frobenius norm error bound obtained in the existing work. 
\end{abstract}


\section{Introduction}\label{sec1}

In the past few decades, most existing work for system identification has been investigated in high-dimensional setting where the ambient dimension is much larger than the sample size. The motivations for these work originates from the fact that the data generated from modern science and engineer field is extremely large, often with the dimension comparable to or possibly much larger than the sample size. In such settings, the classical asymptotic theory often fails to provide useful predictions without imposing some low-dimensional structural constraints. Consequently, there are a lot of work based on different types of constraints. A growing body of literature has focused on the system with sparsity constraints, including sparse linear regression \cite{wai2009,buh2011}, sparse regression matrices \cite{bas2015,per2010,san2013}.

In this paper, we focus on the identification problem for dynamic systems. In general, such problem can be cast as the multivariate linear regression problem \cite{dea2017,gho2018,hal2016,yuan2012}
\begin{equation}\label{int1}
  X=Z\Theta^\ast+W,
\end{equation}
where $X:=[x_1,x_2,\cdots,x_N]^T\in \mathbb{R}^{N\times q}$ and $Z:=[z_1,z_2,\cdots,z_N]^T\in \mathbb{R}^{N\times p}$ consist of the data of response variable $x_i\in \mathbb{R}^{q}$ and explanatory variable $z_i\in \mathbb{R}^{p}$, respectively, $\Theta^\ast \in \mathbb{R}^{p\times q}$ is the coefficient matrix, $W:=[w_1,w_2,\cdots,w_N]^T\in \mathbb{R}^{N\times q}$ is the matrix consisting of noise vectors $w_i$. The multivariate regression problem has been widely applied in macroeconomics \cite{koop2010}, neuroscience \cite{vou2010} and many other areas of applications for predicting response variables. Our interest in this paper is the problem for estimating the coefficient matrix $\Theta^\ast$ that is exactly low-rank (or approximately low-rank). The classical and common method to estimating the coefficient matrix $\Theta^\ast$ is least-squares method \cite{chal2018,dea2017}. However, this method fails to estimate the low-rank matrix $\Theta^\ast$ because its optimal solution is non-unique and will seriously overfit the sample data in high-dimensional setting.

A large number of methods have been proposed to overcome these problems. A natural optimization program is via rank-minimization \cite{mes1998}. Later the selection criterion based on the least-squares cost function with rank penalty is proposed to obtain low-rank solutions \cite{bun2011}. Unfortunately, however, the rank constraint makes it difficult to solve these nonconvex problems. This difficulty motivates replacing the rank penalty with other tractable penalty.
   Accordingly, \cite{faz2002} introduced the nuclear norm heuristics to solve the rank minimization problem. These heuristics are numerically very efficient to obtain low-rank solutions. Furthermore, theoretical properties of these heuristics have been investigated over the past few years. In noiseless setting ($W=0$), there are many studies focused on exact recovery of $\Theta^\ast$. \cite{gro2011} showed that under ``incoherence condition" the low-rank matrix can be recovered exactly from $\Omega \left(r(p+q)\log^2 p\right)$ noisless data ($r$ is the rank of matrix $\Theta^\ast$).  \cite{rec2010} introduced the restricted isometry property (RIP) under which the nuclear norm heuristics can be guaranteed to exactly recover the low-rank solution from $\Omega\left(r(p+q)\log p\right)$ noiseless data. In the noisy setting, many studies focused on the estimation problem of $\Theta^\ast$. \cite{can2009} showed that under RIP, the constrained nuclear norm minimization can stably recover the low rank matrix from $\Omega\left((p+q)r \right)$ noisy data. \cite{neg2009} introduced the weaker and more general restricted strong convexity condition (RSC), and established a consistency and convergence rates for the estimator under certain settings. \cite{neg2011} analyzed the nuclear norm relaxation, showing that the nonasymptotic error bounds on the Frobenius norm hold under RSC. It is noting that much of existing work is based on the RIP or RSC.
    However, depending on the measurement matrix, the RIP or RSC may or may not applicable. For example, the random Pauli measurements suggested in \cite{gro2010} is unknown whether the RIP {holds}\footnote{Throughout the paper, the RIP or RSC holds in the sense of `` high probability ''.}. Moreover, for complex dynamic systems, the rows/columns of corresponding measurement matrix are intrinsically related due to the presence of temporal dependence across observation. This makes it difficult to verify the RIP or RSC of the measurement matrix.

To study the identification problem for dynamic systems, we consider a Vector Auto-Regressive with eXogenous variables (VARX) model. Given the input-output data, the VARX model can be rewritten as Eq.~\eqref{int1} where the measurement matrix is row-related. It is noting that many studies in the past few years have assumed that the measurement matrix had a typical random structures. The common random structures include Gaussian measurement ensemble and row-Gaussian measurement ensemble, all of which are row-independent.   The Gaussian measurement ensemble have been shown to satisfy RIP with high probability, under which the nuclear norm heuristics can recover the low rank solution from noisy measurement \cite{can2009}. Although the row-Gaussian measurement ensemble generally fails to satisfy RIP, it was proved to satisfy the weaker RSC \cite{neg2011}. Note that most of the above studies are based on an observation operator, but here we focus on a observation matrix. Hence, we perform independent repeated sampling for eliminating the dependence of rows of the measurement matrix. This way is important for potentially unstable system. More generally, we assume that the dynamic system noise has a sub-Gaussian tail behaviors \cite{ric2012}. Subsequently, exciting this dynamical systems by sub-Gaussian input, the state of this system at any time has a sub-Gaussian tail behaviors. Consequently, by considering only a state of fixed moment, we can obtain a multivariate linear regression problem \eqref{int1} with the row-sub-Gaussian measurement matrix where its rows are independent isomorphic sub-Gaussian. The row-sub-Gaussian matrix is more general than the Gaussian matrix and the row-Gaussian matrix. However, such matrix generally does not satisfy RIP. In view of this, we introduce the notion of weak RIP, which is weaker than RIP and stronger than RSC, and show that it holds with high probability for the row-sub-Gaussian matrix. Thereby, the nuclear-norm minimization can recover the low-rank exactly with high probability from $\Omega((m+n)r)$ noiseless data (where $r$, $m$ and $n$ represent the rank of the coefficient matrix, the dimensions of system input and output respectively). Further, for noisy observations, we propose a least squares with nuclear norm penalty. Under RSC assumption, such nuclear norm heuristics can derive nonasymptotic bounds on the Frobenius norm error \cite{neg2011}. We now turn to alternative form of restricted curvature \cite{wai2015}, which involves the gradient of the cost function. Based on this condition, we can obtain a tighter error bound. In particular, we introduce the notion of the operator-norm curvature condition. Later, we show that the operator-norm curvature condition holds with high probability for the row-sub-Gaussian matrix. Thereout, specifying a suitable choice of the regularization parameter, we obtain an operator norm error bound $\mathcal{O}(\sqrt{(m+n)/N})$ (where $N$ represents the number of samples) from $\Omega(m+n)$ noisy datasets with high probability.

This rest of this paper is organized as follows. We begin in Section \ref{sec2} with the model setup and the corresponding estimation procedure. We then in Section \ref{sec3} devote to the statements of our theoretical results on the asymptotic behaviors of nuclear norm heuristic methods, which is established in two different scenarios. Finally, we conclude our paper in Section \ref{conclusion}. Detailed technical proofs are provided in the Appendix.

\textbf{Notations.} Given a matrix $A\in \mathbb{R}^{d_1\times d_2}$, we let $\mathrm{rowspan(A)}\subseteq \mathbb{R}^{d_2}$ and $\mathrm{colspan(A)}\subseteq \mathbb{R}^{d_1}$ be its row and column spaces respectively. We denote by  $A_{i,:}$ for any $1\leq i\leq d_1$ the $i^{th}$ row of $A$ and write $A_{:,j}$ for its $j^{th}$ column, $1\leq j\leq d_2$.  We denote by $\sigma_1(A) \geq \sigma_2(A) \geq \cdots \geq 0$ the singular values of $A$. Based on these singular values, we define the matrix norms, including the nuclear norm $\|A\|_{\mathrm{nuc}}=\sum \limits_{i=1}^{d}\sigma_i(A)$, the operator norm $\|A\|_{\mathrm{op}}=\sigma_1(A)$ and the Frobenius norm $\|A\|_{\mathrm{\mathrm{F}}}=\sqrt{\mathrm{trace}(A^TA)}=\sqrt{\sum \limits_{i=1}^{d}\sigma_i^2(A)}$, where $d:=\min\{d_1,d_2\}$ and the superscript ``T" denotes the matrix transpose.
 For a matrix $\Delta\in \mathbb{R}^{d_1\times d_2} $ and a subspace $\mathcal{M}\subseteq \mathbb{R}^{d_1\times d_2}$, we use the notation $\Delta_{\mathcal{M}}$ to refer to the projection of the matrix $\Delta$ onto the subspace $\mathcal{M}$. Given a vector $u\in \mathbb{R}^d$, we let $\|u\|_2$ denote the Euclidean norm (or $l_2$-norm) of $u$. We denote by $\langle \cdot,\cdot \rangle$ the inner product with respect to $\ell_2$-norm. For any integer $d\geq 1$, we denote by $\mathbb{S}^{d-1}$  the unit sphere with respect to $\ell_2$-norm in $\mathbb{R}^d$. The notation $f_1(n)=\Omega(f_2(n))$ represents that there exists a constant $c$ such that $f_1(n)\geq cf_2(n)$, while the notation $f_1(n)=\mathcal{O}(f_2(n))$ represents that there exists a constant $c$ such that $f_1(n)\leq cf_2(n)$. For a set $S$, we denote by $|S|$ the cardinality of $S$.

\section{Problem formulation}\label{sec2}

In this section, we consider a system with linear state-space form
\begin{equation}\label{p1}
x(t+1)=Ax(t)+Bu(t)+w(t),
\end{equation}
where $x(t)\in \mathbb{R}^n$, $u(t)\in \mathbb{R}^m$, $w(t)\in \mathbb{R}^{n}$ are the system state, input, noise at time instant $t$, respectively, $A\in \mathbb{R}^{n\times n}$ is the state transition matrix, $B\in \mathbb{R}^{n\times m}$ is the input matrix. Our objective is to estimate the coefficient matrices under high-dimensional scaling. However, it is impossible to estimate the coefficient matrices in such setting unless the model is equipped with some low-dimensional structures. In this work, we assume that the matrices $A$, $B$ (or block matrix $[A,B]$) are low-rank.
Note that, the VAR model, as a special case, has been studied \cite{lut2005,lin2017}. We then estimate the coefficient matrices by a finite number of sample sets like $\big{\{}(x(t),u(t))\big{\}}_{t=0}^{T^{0}}$, where $T^0$ represents the maximal sampling time. To this end, we repeat the sampling scheme over the interval $[0,T^0]$ for obtaining multiple sample sets.

For a sample set  $\big{\{}(x(t),u(t))\big{\}}_{t=0}^{T^{0}}$, let $z(t):=\left[
                                                                                                                                                                                                                             \begin{array}{c}
                                                                                                                                                                                                                               x(t) \\
                                                                                                                                                                                                                               u(t) \\
                                                                                                                                                                                                                             \end{array}
                                                                                                                                                                                                                           \right]\in \mathbb{R}^{n+m}$,
 $\Theta^\ast:=[A,B]^T\in \mathbb{R}^{(n+m)\times n}$, one has
\begin{equation}\label{p5}
x^T (t+1)=z^T(t)\Theta^\ast+w^T(t).
\end{equation}
Stacking \eqref{p5} for $t=0,1,\cdots,T^0-1$, one has
\begin{equation}\label{p2}
  X_{T^0}=Z_{T^0}\Theta^\ast+W_{T^0},
\end{equation}
where $X_{T^0}:=[x(1),x(2),\cdots,x(T^0)]^T, Z_{T^0}:=[z(0),z(1),\cdots,z(T^0-1)]^T$, and $ W_{T^0}:=[w(0),w(1),\cdots,w(T^0-1)]^T$.

 We excite the system by arranging a system input to obtain $N$ sample sets with length $T^0$. Note though that the collection process of each sample set is independent. To make a distinguish, we denote by $\{(x^{(i)}(t),u^{(i)}(t))\}_{t=0}^{T^0}$ the $i^{th}$ sample set. Correspondingly, the representation \eqref{p2} for the $i^{th}$ sample set is written as $X^{(i)}_{T^0}=Z^{(i)}_{T^0}\Theta^\ast+W^{(i)}_{T^0}$. For lightening the notation, we omit the subscript $T_0$. Based on these sample sets, the commonly adopted method to identify the coefficient matrix is the least squares method (LS), through which we obtain the least squares estimator
\begin{equation}\label{p3}
  \widehat{\Theta}_{LS}=\arg \min \limits_{\Theta} \frac {1}{N} \sum \limits_{i=1}^N\|X^{(i)}-Z^{(i)}\Theta\|_F ^2.
\end{equation}
If the sum $\sum \limits_{i=1}^N  (Z^{(i)})^TZ^{(i)}$ is invertible, the optimal solution of \eqref{p3} is unique, which can be computed as
\[
 \widehat{\Theta}_{LS}=\left( \sum \limits_{i=1}^N  (Z^{(i)})^TZ^{(i)} \right)^{-1}\left( \sum \limits_{i=1}^N  (Z^{(i)})^T X^{(i)} \right),
\]
Then the error matrix is
\[
 \widehat{\Theta}_{LS}-\Theta^\ast=\left( \sum \limits_{i=1}^N  (Z^{(i)})^TZ^{(i)} \right)^{-1}\left( \sum \limits_{i=1}^N  (Z^{(i)})^T W^{(i)} \right).
\]
Specially, the coefficient matrix can be exactly recovered by applying LS in the absence of noise ($W$= 0). Note that the state $x^{(i)}(t)$ at time $t$ is influenced by the noise at time $0,1,\cdots, t-1$, and hence $X^{(i)}$ and $W^{(i)}$ are related. It leads to the difficulty of studying and analyzing the behavior of estimators $\widehat{\Theta}_{LS}$. To alleviate this problem, we consider the state at time $T^0$ instead of the whole sample set. Further, we define
$X:=[ x^{(1)}(T^0), x^{(2)}(T^0),\cdots, x^{(N)}(T^0)]^T\in \mathbb{R}^{N\times n}$, $Z:=[ z^{(1)}(T^0-1), z^{(2)}(T^0-1),\cdots, z^{(N)}(T^0-1)]^T \in \mathbb{R}^{N\times (n+m)}$, $W:=[ w^{(1)}(T^0-1), w^{(2)}(T^0-1),\cdots, w^{(N)}(T^0-1)]^T\in \mathbb{R}^{(n+m)\times n}$.
Then, one has
\begin{equation}\label{p4}
X=Z\Theta^\ast+W,
\end{equation}
Since the optimal solution of the LS is non-unique in high-dimensional setting, the LS fails to recover the low-rank matrix. To reduce the dimensionality of the model, we use a low-dimensional constraint. Thus, we consider the following convex relaxation program
\begin{equation}\label{m1}
  \widehat{\Theta}=\arg \min \limits_{\Theta} \big{\{}\mathcal{L}_N(\Theta)+\lambda_N\mathcal{R}(\Theta)\big{\}},
\end{equation}
where $\mathcal{L}_N(\Theta):=\frac {1}{2N}\|X-Z\Theta\|_{\mathrm{F}} ^2$, $\mathcal{R}(\cdot):=\|\cdot\|_{\mathrm{nuc}}$ stand for the cost function and nuclear norm, respectively, $\lambda_N$ is a user-defined regularization penalty to balance the strength of the loss and regularization term. The identification problem focuses on determining how many sample sets are needed to recover the low-rank matrix. The theoretical properties of this program will be analyzed in the next section.
\section{Theoretical results}\label{sec3}
 Assume that the coefficient matrix $\Theta^\ast$ has rank $r< n$. Then $\Theta^\ast$ can be decomposed as $\Theta^\ast =UDV^T$, where the diagonal matrix $D\in \mathbb{R}^{r\times r}$ has the $r$ nonzero singular values of $\Theta^\ast$ in its diagonal entries, $U\in \mathbb{R}^{(n+m)\times r}$ and $V\in \mathbb{R}^{n\times r}$ are orthogonal matrices, with their columns corresponding to the left and right singular vectors of $\Theta^\ast$, respectively. Define two subspaces $\mathcal{M}$, $\mathcal{\bar{M}}^\perp$ of $\mathbb{R}^{(n+m)\times n}$ with
\[
\mathrm{\mathcal{M}=\big{\{}\Theta \in \mathbb{R}^{(n+m)\times n}|rowspan(\Theta)\subseteq \mathcal{V},colspan(\Theta)\subseteq \mathcal{U} \big{\}},}
\]
\[
\mathrm{\mathcal{\bar{M}}^\perp=\big{\{}\Theta \in \mathbb{R}^{(n+m)\times n}|rowspan(\Theta)\perp \mathcal{V}, colspan(\Theta)\perp \mathcal{U}\big{\}}},
\]
where $\mathcal{U}$ and $\mathcal{V}$ represent the subspaces generated by the columns of matrices $U$ and $V$ respectively. Here,
$\mathrm{\mathcal{\bar{M}}}$ denotes the subspace orthogonal to $\mathrm{\mathcal{\bar{M}}^\perp}$.
Obviously, $\Theta^\ast \in \mathcal{M}\subseteq \mathrm{\mathcal{\bar{M}}}$.

This section will make use of the standard results on sub-Gaussian random vector and random matrix theory for obtaining probabilistic statements.
\begin{definition}\label{def2}
(Sub-Gaussian vector).\ A random vector $x\in \mathbb{R}^d$ with zero mean is \emph{sub-Gaussian} with parameter $\sigma>0$ if for any fixed $v\in \mathbb{S}^{d-1}$, $$\mathbb{E}e^{\lambda \langle v,x \rangle} \leq e^{\frac {\lambda^2 \sigma^2}{2}}\ \ \ \  \mathrm{for}\ \mathrm{all}\  \lambda \in \mathbb{R}.$$
In particular, when $d=1$, we call it \emph{sub-Gaussian variable}.
\end{definition}
\begin{remark}\label{remark1}
\emph{\emph{}There are several examples considered frequently for sub-Gaussian vector.}
\emph{\begin{itemize}
 \item Random vector $x\in \mathbb{R}^d\thicksim \mathcal{N}(0,\Sigma)$. This is a more common situation. Since $\langle v,x\rangle\thicksim \mathcal{N}(0,v^T\Sigma v)$ and $v^T\Sigma v\leq \|\Sigma\|_{\mathrm{op}}$ for each $v\in \mathbb{S}^{d-1}$, the vector $x$ is sub-Gaussian with parameter at most $\sigma^2=\|\Sigma\|_{\mathrm{op}}$.
  \item Random vector $x\in \mathbb{R}^d$ has independently identically distributed (i.i.d.) entries, where each entry $x_i$ is zero-mean and sub-Gaussian. For example, the Gaussian distribution ($x_i\thicksim\mathcal{N}(0,\sigma^2)$);
   the uniform distribution on the interval $[-a,a]$. For all of these cases, the random variable $\langle v,x\rangle$ is sub-Gaussian variable with any vector $v\in \mathbb{S}^{d-1}$.
\end{itemize}}
 \end{remark}
 In order to establish a probability framework, we initialize this system by $x(0)=0$ and assume the noise (if exists) has a sub-Gaussian tail behaviors. Using sub-Gaussian excitation, then the measure matrix $Z$ is row-sub-Gaussian (see Propositions \ref{pro2} and \ref{pro1}). The row-sub-Gaussian matrices are a more general random matrices that its rows are i.i.d. sub-Gaussian vectors. Such matrices include Gaussian matrix \cite{can2013}, Bernoulli matrix \cite{men2008}, and more generally the matrix with sub-Gaussian entries \cite{ver2000}.
 We then will analyze the nonasymptotic behaviors of the estimator of the low-rank matrix $\Theta^\ast$. We begin with several lemmas.
  \begin{lemma}\label{lem7}\cite{rec2010}
  Nuclear norm $\|\cdot\|_{\mathrm{nuc}}$ is decomposable with respect to $\left(\mathcal{M},\bar{\mathcal{M}}^\perp\right)$, that is,
  \[
  \|A+B\|_{\mathrm{nuc}}=\|A\|_{\mathrm{nuc}}+\|B\|_{\mathrm{nuc}}
  \]
 holds for all $A\in \mathcal{M}$, $B\in \bar{\mathcal{M}}^\perp$.
  \end{lemma}

\begin{lemma}\label{lem1}\cite{wai2015}
Let $Z\in \mathbb{R}^{N\times (n+m)}$ be a row-sub-Gaussian matrix and $\Sigma$ be its covariance matrix. Then there exists a constant $\beta>0$ such that the sample covariance matrix $\widehat{\Sigma}:=\frac {1}{N}Z^T Z$ satisfies the bounds
\begin{equation}\notag
  \mathbb{P}\left[\|\widehat{\Sigma}-\Sigma\|_{\emph{op}} \geq 16\sqrt{6}\beta^2 \left( \sqrt {\frac {n+m}{N}}+\frac {n+m}{N}\right)+\delta\beta^2\right]\leq e^{-N \min
  \big{\{}\frac{\delta}{16\sqrt{2}},\frac{\delta^2}{512}\big{\}}}.\label{m3}
\end{equation}
\end{lemma}

\begin{remark}
In particular, for the standard Gaussian matrix $Z$, one has
\[
\|\widehat{\Sigma}-I_{n+m}\|_{\mathrm{op}}  \preceq \sqrt {\frac {n+m}{N}}+\frac {n+m}{N}
\]
with high probability. Similar results can be seen in \cite{ver2010,wai2009}. 
\end{remark}

  In the identification process for the dynamic system \eqref{p1}, if $u^{(i)}(t)$ and $w^{(i)}(t)$ are independent Gaussian vectors, then $z^{(i)}(T^0-1)\in \mathbb{R}^{n+m}$ is a Gaussian vector, and hence the above concentration inequality holds for the measurement matrix. In the next section, we will apply the above concentration inequality to obtain the results for the recovery of the low-rank matrix using nuclear norm heuristic method.

 \subsection{Weak Restricted Isometry and Recovery of Low-rank Matrices}
    We begin by focusing on the scenario where the observations are perfect or noiseless. In such setting, our aim is to find a low-rank matrix $\Theta$ such that $X=Z\Theta$. Obviously, this problem can be cast as the optimization problem
\begin{equation}\label{noise1}
  \min \limits_{\Theta} \mathrm{rank}(\Theta)\ \ \ \ \ \ \mathrm{such}\  \mathrm{that}\ X=Z\Theta,
\end{equation}
where the definitions of $X$, $Z$ are shown in \eqref{p4}. It is known to be NP-hard. To alleviate this problem, we replace the rank constraints by nuclear norm constraints, which leading to the convex program
 \begin{equation}\label{noise5}
     \min \limits_{\Theta} \|\Theta\|_{\mathrm{nuc}}\ \ \ \ \ \  \mathrm{such}\  \mathrm{that}\ X=Z\Theta.
 \end{equation}
 Let $\widehat{\Theta}$ be an optimal solution to the program \eqref{noise5}. In this section, we will characterize specific cases when we can a priori guarantee that $\widehat{\Theta}=\Theta^\ast$. Since the row-sub-Gaussian matrix does not satisfy the RIP, we introduce the weak RIP that holds with high probability for the row-sub-Gaussian matrix.

\begin{definition}\label{def3}
For each integer $r$ with $r< n$, we say that the matrix $Z\in \mathbb{R}^{N\times (n+m)}$ satisfies \emph{weak RIP} of order $r$ with constants $\delta\in (0,1)$, $K_2\geq K_1>0$ (does not depends on $r$) if
\begin{equation}\label{noise4}
  K_1(1-\delta)\|\Delta\|_{\mathrm{F}} \leq \frac {\|Z\Delta\|_{\mathrm{F}}}{\sqrt{N}} \leq K_2(1+\delta)\|\Delta\|_{\mathrm{F}}
\end{equation}
holds for all matrices $\Delta\in \mathbb{R}^{(n+m)\times n}$ of rank at most $r$. We let $\delta_r(K_1,K_2)$ denote the the weak-RIP constant of order $r$, i.e., the smallest number $\delta>0$ such that \eqref{noise4} holds.
\end{definition}
\begin{remark}
The constants $K_1$, $K_2$ in \eqref{noise4} only depend the matrix $Z$. Then, by the definition of weak RIP, one has $\delta_r(K_1,K_2)\leq  \delta_{r^{\prime}}(K_1,K_2)$ for $r\leq r^{\prime}$. In particular, when $K_1=K_2$, the weak RIP is the common restricted isometry property. Here, we omit the case where $K_1>K_2$ because it can be converted into the case where $K_1=K_2$. In addition, the weak RIP is stronger than the RSC due to the existence of the non-isometry upper bounds.
\end{remark}
Define
\begin{equation*}
s:=
\begin{cases}
1,& \text{if $K_1=K_2,$}\\
\textbf{[}(\frac{K_2}{K_1})^2\textbf{]}+1,& \text{if $K_1<K_2$},
\end{cases}
\end{equation*}
where $\textbf{[}\cdot \textbf{]}  $ stands for the rounding function that its function value in $x$ is the maximum integer that does not exceed $x$. Obviously, $s\geq (K_2/K_1)^2\geq 1$. Based on the analysis of the weak RIP, we obtain the following two recovery conclusions.

\begin{proposition}\label{pro3}
If the weak-RIP constant of order $2r$ satisfies $\delta_{2r}(K_1,K_2)<1$, then $\Theta^\ast$ is the only matrix of rank at most $r$ satisfying $Z\Theta=X$.
\end{proposition}
\noindent {\bf Proof.}\ Obviously, $\Theta^\ast$ is a matrix of rank $r$ satisfying $Z\Theta=X$. Then, we only need to prove its uniqueness. Using the proof by contradiction, we assume that there exists a matrix $\Theta_0$ of rank $r$ satisfying $Z\Theta_0=X$ and $\Theta_0\neq \Theta^\ast$. Let $\Delta:=\Theta_0-\Theta^\ast$. Then, $\Delta$ has rank at most $2r$ and $Z\Delta=0$. By the weak RIP, $0=\frac{\|Z\Delta\|_{\mathrm{F}}^2}{N}\geq K_1\left(1-\delta_{2r}(K_1,K_2)\right)\|\Delta\|_{\mathrm{F}}^2$, which implies $\Delta=0$, i.e. $\Theta_0= \Theta^\ast$. This is a contradiction. \qed

\begin{theorem}\label{th5}
For the integer $r\geq 1$, we have $\widehat{\Theta}=\Theta^\ast$ if $\delta_{(2+3s)r}(K_1,K_2)<5-2\sqrt{6}$.
\end{theorem}
\noindent {\bf Proof.}\ By the definition of $\widehat{\Theta}$, one has $\|\widehat{\Theta}\|_{\mathrm{nuc}}\leq \|\Theta^\ast\|_{\mathrm{nuc}} \leq \|\Theta^\ast_{\mathcal{M}}\|_{\mathrm{nuc}}$. Let $\widehat{\Delta}:=\widehat{\Theta}-\Theta^\ast$. Then
\begin{align*}
  \|\Theta^\ast\|_{\mathrm{nuc}} & \geq   \|\Theta^\ast+\widehat{\Delta}\|_{\mathrm{nuc}}= \|\Theta^\ast_{\mathcal{M}}+\widehat{\Delta}_{\bar{\mathcal{M}}}+\widehat{\Delta}_{\bar{\mathcal{M}}^\perp}\|_{\mathrm{nuc}}    \\
   & \overset {(i)}{\geq} \|\Theta^\ast_{\mathcal{M}}+\widehat{\Delta}_{\bar{\mathcal{M}}^\perp}\|_{\mathrm{nuc}}-\|\widehat{\Delta}_{\bar{\mathcal{M}}}\|_{\mathrm{nuc}}\\
   & \overset {(ii)}{=} \|\Theta^\ast_{\mathcal{M}}\|_{\mathrm{nuc}}+\|\widehat{\Delta}_{\bar{\mathcal{M}}^\perp}\|_{\mathrm{nuc}}-\|\widehat{\Delta}_{\bar{\mathcal{M}}}\|_{\mathrm{nuc}},
\end{align*}
where inequality $(i)$ follows from the triangle inequality and equality $(ii)$ follows from Lemma \ref{lem7}. Thereby, we conclude that
\begin{equation}\label{noise6}
\|\widehat{\Delta}_{\bar{\mathcal{M}}^\perp}\|_{\mathrm{nuc}}\leq \|\widehat{\Delta}_{\bar{\mathcal{M}}}\|_{\mathrm{nuc}}.
\end{equation}
We write the Singular Value Decomposition (SVD) as $\widehat{\Delta}_{\bar{\mathcal{M}}^\perp}:=UD V^T$, where $U:=[u_1,u_2,\cdots,u_n]\in \mathbb{R}^{(n+m)\times n}$ and $V:=[v_1,v_2,\cdots,v_n]\in \mathbb{R}^{n\times n}$ are the column orthogonal matrices, $D:=\mathrm{diag}(\sigma_1,\sigma_2,\cdots,\sigma_n) $ is a diagonal matrix that its diagonal entries are the singular values of the matrix $\widehat{\Delta}_{\bar{\mathcal{M}}^\perp}$ in decreasing order. Let $n=(3sr)p+q$, where $p,q$ are positive integers and $0\leq q< 3sr$. Define the index set $I_i=\{(3sr)(i-1)+1,\cdots,(3sr)i\}$ for any $i=1,\cdots,p$ and $I_{p+1}=\{(3sr)p+1,\cdots,n\}$. Then the matrix $\widehat{\Delta}_{\bar{\mathcal{M}}^\perp}$ can be decomposed as $\widehat{\Delta}_{\bar{\mathcal{M}}^\perp}:=\sum \limits_{i=1}^{p+1}\Delta_i$, where $\Delta_i:=\sum \limits_{j\in I_i }\sigma_i u_iv_i^T$ and its rank is at most $3sr$. Obviously, $\Delta_i^T\Delta_j=0$ for any $i\neq j$. Assume that $\sigma_1\geq \sigma_2\geq \cdots\geq \sigma_n\geq 0$. Then, one has
\[
\sigma_k\leq \frac {1}{3sr}\sum \limits_{j\in I_i}\sigma_j=\frac {1}{3sr}\|\Delta_i\|_{\mathrm{nuc}}
\]
for any $k\in I_{i+1}$. Hence,
\[
\|\Delta_{i+1}\|_{\mathrm{F}}^2=\sum \limits_{k\in I_{i+1}} \sigma_k^2\leq 3sr\max \limits_{k\in I_{i+1}}\sigma_k^2\leq  \frac {1}{3sr}\|\Delta_i\|_{\mathrm{nuc}}^2.
\]
Thereafter,
\begin{equation*}
  \sum \limits_{i=2}^{p+1}\|\Delta_{i}\|_{\mathrm{F}}\leq \frac {1}{\sqrt{3sr}}\sum \limits_{i=1}^{p+1}\|\Delta_{i}\|_{\mathrm{nuc}}=\frac {1}{\sqrt{3sr}}\|\widehat{\Delta}_{\bar{\mathcal{M}}^\perp}\|_{\mathrm{nuc}}\leq \frac {\sqrt{2r}}{\sqrt{3sr}}\|\widehat{\Delta}_{\bar{\mathcal{M}}^\perp}\|_{\mathrm{F}},
\end{equation*}
where last inequality follows from the fact the $\mathrm{rank}(\widehat{\Delta}_{\bar{\mathcal{M}}})\leq 2r$ (any matrix in $\bar{\mathcal{M}}$ has rank at most 2$r$). Since the rank of $\widehat{\Delta}_{\bar{\mathcal{M}}}+\Delta_1$ is at most $(2+3s)r$, we have
\begin{align*}
\frac {\|Z\widehat{\Delta}\|_{\mathrm{F}}}{\sqrt{N}} \geq & \frac {\|Z(\widehat{\Delta}_{\bar{\mathcal{M}}}+\Delta_1)\|_{\mathrm{F}}}{\sqrt{N}}-\frac {\sum \limits_{i=2}^{p+1}\|Z\Delta_i\|_{\mathrm{F}}}{\sqrt{N}}       \\
   \geq & K_1\left(1-\delta_{(2+3s)r}(K_1,K_2)\right)\|\widehat{\Delta}_{\bar{\mathcal{M}}}+\Delta_1\|_{\mathrm{F}} -K_2\left(1+\delta_{3sr}(K_1,K_2)\right)\sum \limits_{i=2}^{p+1}\|\Delta_i\|_{\mathrm{F}}\\
\geq &\left[ K_1\left(1-\delta_{(2+3s)r}(K_1,K_2)\right)-K_2\frac {\sqrt{2r}}{\sqrt{3sr}}\left(1+\delta_{3sr}(K_1,K_2)\right)\right]\|\widehat{\Delta}_{\bar{\mathcal{M}}}\|_{\mathrm{F}}\\
\geq & K_1\left[\left(1-\delta_{(2+3s)r}(K_1,K_2)\right)- \sqrt{\frac {2}{3}}\left(1+\delta_{3sr}(K_1,K_2)\right)  \right]\|\widehat{\Delta}_{\bar{\mathcal{M}}}\|_{\mathrm{F}}.
\end{align*}
Observe that $Z\widehat{\Delta}=Z\widehat{\Theta}-Z\Theta^\ast=0$. Then $\widehat{\Delta}_{\bar{\mathcal{M}}}=0$ when $\left(1-\sqrt{\frac {2}{3}}\right)-\left(\delta_{(2+3s)r}(K_1,K_2)+\sqrt{\frac {2}{3}}\delta_{3sr}(K_1,K_2)\right)>0$.
 Finally, applying the monotonicity for weak-RIP constant yields the desired result.\qed

The above conclusions provide the conditions that guarantee $\widehat{\Theta}=\Theta^\ast$. In the following we will demonstrate that the measurement matrix $Z$ satisfies these conditions with overwhelming probability.

By the state equation \eqref{p1}, one has
\begin{equation}\label{noise2}
  x^{(i)}(T^0-1)= A^{T^0-2}Bu^{(i)}(0)+\cdots+Bu^{(i)}(T^0-2).
 \end{equation}
 We take into account the special excitation for facilitating the analysis of the statistical properties of estimators $\widehat{\Theta}$. Using the sub-Gaussian excitation with parameter $\sigma_u$, we then have the following results.

 \begin{proposition}\label{pro2}
 The random vector $z^{(i)}(T^0-1)\in \mathbb{R}^{n+m}$ is sub-Gaussian with parameter $\sigma_z$ for any $i=1,2,\cdots,N$, where $\sigma_z=\sqrt{\sum \limits_{k=0}^{T^0-2}\left(\|A^{T^0-2-k}B\|_{\emph{op}}\sigma_u^2\right)+\sigma_u^2}$.
 \end{proposition}
\noindent {\bf Proof.}\ For any fixed $v\in \mathbb{R}^{n+m}$ and any $\lambda \in \mathbb{R}$, by Eq. \eqref{noise2} and the independence of $u^{(i)}(i=1,\cdots,T^0-1)$, we have
\begin{align*}
  \mathbb{E}[e^{\lambda \langle v,z^{(i)}(T^0-1)\rangle}]=&\mathbb{E}e^{\lambda \langle v_1,x^{(i)}(T^0-1) \rangle}\mathbb{E} e^{\lambda \langle v_2,u^{(i)}(T^0-1) \rangle}\\
  =& \mathop {\prod} \limits_{k=1}^{T^0-2}\left(\mathbb{E}e^{\lambda \langle v_1,A^{T^0-2-k}Bu^{(i)}(k) \rangle}\right)\mathbb{E} e^{\lambda \langle v_2,u^{(i)}(T^0-1) \rangle}\\
   =& \mathop {\prod} \limits_{k=1}^{T^0-2} \left(\mathbb{E}e^{\lambda \langle (A^{T^0-2-k}B)^Tv_1,u^{(i)}(k) \rangle}\right)\mathbb{E}e^{\lambda \langle v_2,u^{(i)}(T^0-1) \rangle}\\
  \leq & e^{\frac{\lambda^2\sigma_z^2}{2}},
\end{align*}
where $v=\left[
           \begin{array}{c}
             v_1 \\
             v_2 \\
           \end{array}
         \right]
$, $v_1\in \mathbb{R}^n$, $v_2\in \mathbb{R}^m$, which shows that $z^{(i)}(T^0-1)$ is sub-Gaussian with parameter $\sigma_z$.    \qed

Denote by $\Sigma$ the covariance matrix of random vector $z^{(i)}(T^0-1)$. Let $\gamma_{\min}(\Sigma)$ and $\gamma_{\max}(\Sigma)$ be the minimal and maximal singular value of matrix $\Sigma$ respectively. Assume that $\Sigma$ is invertible, i.e., $\gamma_{\min}(\Sigma)>0$. Further, without loss of generality, we assume that $\gamma_{\max}(\Sigma)\leq 1$. Indeed, if the $\gamma_{\max}(\Sigma)> 1$, by standardizing the matrix $\Sigma$, i.e.,  $\hat{Z}:=Z/\sqrt{\gamma_{\max}(\Sigma)}$, it can come down to the case where $\gamma_{\max}(\Sigma)\leq 1$.

 In \cite{rec2010}, the nearly isometric random matrices are proven to obey the RIP. A random matrix $Y\in \mathbb{R}^{N \times (n+m)}$ to be nearly isometric has to satisfy two conditions. First, it is isometric in expectation, i.e., $\mathbb{E}\|Y\Theta\|^2_\mathrm{F}=\|\Theta\|^2_\mathrm{F}$ for any matrices $\Theta\in \mathbb{R}^{(n+m)\times n}$. There are some matrices that satisfy this property. For example, the Gaussian measurement matrix with i.i.d. $\mathcal{N}(0,1/(m+n))$ entries. Second, the probability of large deviation of the length is exponentially small. However, the row-sub-Gaussian matrix generally does not satisfies these conditions. Here we will demonstrate that the row-sub-Gaussian matrix satisfies a more general deviation inequality.
\begin{lemma}\label{lem6}
For any fixed $t>0$, there exist constants $c_1$, $c_2$ such that when $N\geq c_1(n+m)$,
\[
 \mathbb{P}\left[ \frac {1}{N}\left|\|Z\Theta\|_{\mathrm{F}}^2-\mathbb{E}\|Z\Theta\|_{\mathrm{F}}^2\right|\geq  t \|\Theta\|_{\mathrm{F}}^2\right]\leq n \exp(-c_2N).
 \]
 \end{lemma}
 \noindent {\bf Proof.}\  Observe that
\begin{align*}
    \mathbb{P}\left[ \frac {1}{N}\left|\|Z\Theta\|_{\mathrm{F}}^2-\mathbb{E}\|Z\Theta\|_{\mathrm{F}}^2\right|\geq  t \|\Theta\|_{\mathrm{F}}^2\right]\leq & \mathbb{P}\left[\sum \limits_{k=1}^{n}\left|\Theta_{:,k}^T(\widehat{\Sigma}-\Sigma)\Theta_{:,k}\right|\geq  t \|\Theta_{:,k}\|_2^2\right] \\
  \leq & \sum \limits_{k\in \Omega}\mathbb{P}\left[|\Theta_{:,k}^T(\widehat{\Sigma}-\Sigma)\Theta_{:,k}|\geq  t \|\Theta_{:,k}\|_2^2\right] \\
  \leq & n\mathbb{P}\left[\|\widehat{\Sigma}-\Sigma\|_{\mathrm{op}}\geq t\right],
\end{align*}
where $\widehat{\Sigma}:=\frac {1}{N}Z^TZ$ represents the sample covariance matrix, $\Omega:=\{k\in \{1,2,\cdots,n\}\ |\ \Theta_{:,k}\neq 0\}$ is a set composed of non-zero column indices of $\Theta$. By Lemma \ref{lem1}, there exists constants $c_1$, $c_2$ such that $$\mathbb{P}\left[\|\widehat{\Sigma}-\Sigma\|_{\mathrm{op}}\geq t\right]\leq \exp(-c_2N),$$
when $N\geq c_1(n+m)$. Hence, we derive the desired claim.\qed

Using the above technical lemmas and theorems, we can then obtain the following results.
\begin{theorem} \label{th4}
For any fixed $0\leq\delta<1$, there are constants $c$, $d$ such that if $N\geq c(n+m)r$, $Z$ satisfies the weak RIP condition with constants $\delta_r(K_1,K_2)\leq \delta$, $K_1=\sqrt{\gamma_{\min }(\Sigma)}$ and $K_2=\sqrt{\gamma_{\max}(\Sigma)}$ with probability at least $1-e^{-dN}$.
\end{theorem}
 \noindent {\bf Proof.}\ To prove Theorem \ref{th4}, it suffices to show that
 \begin{equation}\label{m15}
 \sqrt{\gamma_{\min}(\Sigma)}(1-\delta) \leq \frac {\|Z\Delta\|_{\mathrm{F}}}{\sqrt{N}} \leq \sqrt{\gamma_{\max}(\Sigma)}(1+\delta)
\end{equation}
holds for all matrices $\Delta\in \mathcal{U}_r:=\{\Delta\in \mathbb{R}^{(n+m)\times n}|\ \mathrm{rank}(\Delta)\leq r,\ \|\Delta\|_{\mathrm{F}}=1 \}$.
For any $\epsilon>0$, by the covering number theorem \cite{can2009} for the set of low-rank matrices, there exists an $\epsilon$-net $\mathcal{V}_r\subseteq \mathcal{U}_r$ with respect to the Frobenius norm, which contains at most $(9/\epsilon)^{(2n+m+1)r}$ elements. Taking $\epsilon=\delta\sqrt{\gamma_{\min}(\Sigma)}/(4\sqrt{2})$, $\mathcal{V}_r$ has at most $(36\sqrt{2}/(\delta\sqrt{\gamma_{\min}(\Sigma)}))^{(2n+m+1)r}$ elements. By Lemma \ref{lem6} with $t=\delta\gamma_{\min}(\Sigma)/2$, we have
\begin{align*}
   & \mathbb{P}\left[\sup \limits_{\Theta\in \mathcal{V}_r}  \frac {1}{N}|\|Z\Theta\|_{\mathrm{F}}^2-\mathbb{E}\|Z\Theta\|_{\mathrm{F}}^2|\geq  \delta\gamma_{\min}(\Sigma)/2 \right] \\
   \leq  & (36\sqrt{2}/(\delta\sqrt{\gamma_{\min}(\Sigma)}))^{(2n+m+1)r}n e^{-c_2N}\\
   \leq & \exp\left(-c_2N+(n+m)r\left(2\log(36\sqrt{2}/(\delta\sqrt{\gamma_{\min}(\Sigma)}))+1\right) \right)\\
   \leq & \exp(-dN),
\end{align*}
when $N\geq c(n+m)r$, where $c:=\max \big{\{}c_1,1+2\log(36\sqrt{2}/(\delta\sqrt{\gamma_{\min}(\Sigma)}))\big{\}}$, $d:=c_2/2$.  Assume that
\begin{equation}\label{m16}
\sup \limits_{\Theta\in \mathcal{V}_r } \frac {1}{N}|\|Z\Theta\|_{\mathrm{F}}^2-\mathbb{E}\|Z\Theta\|_{\mathrm{F}}^2|\leq  \delta\gamma_{\min}(\Sigma)/2.
\end{equation}
 Let
 \[
J:=\sup \limits_{\Theta\in \mathcal{U}_r}\frac {\|Z\Theta\|_{\mathrm{F}}}{\sqrt{N}}.
 \]
  Then there exists a matrix $P\in \mathcal{V}_r $ such that $\|P-\Theta\|_{\mathrm{F}}\leq \delta\sqrt{\gamma_{\min}(\Sigma)}/4\sqrt{2}$. Hence,
 \begin{equation*}
    \frac {1}{\sqrt{N}}\|Z\Theta\|_{\mathrm{F}}\leq \frac {1}{\sqrt{N}}\|Z(P-\Theta)\|_{\mathrm{F}}+\frac {1}{\sqrt{N}}\|ZP\|_{\mathrm{F}} \leq \frac {1}{\sqrt{N}}\|Z(P-\Theta)\|_{\mathrm{F}}+ \sqrt{\gamma_{\max}(\Sigma)}(1+\frac {\delta}{2}),
 \end{equation*}
 where the last inequality follows from the inequality \eqref{m16}.
Let $\Delta:=P-\Theta$. Then $\mathrm{rank}(\Delta)\leq 2r$. By SVD, we can easily get the decomposition : $\Delta=\Delta_1+\Delta_2$, where $\mathrm{rank}(\Delta_i)\leq r$, $i=1,2$ and $\mathrm{tr}((\Delta_1)^T\Delta_2)=0$. Thereby,
\begin{align*}
   \frac {1}{\sqrt{N}}\|Z\Delta\|_{\mathrm{F}}\leq \frac {1}{\sqrt{N}}\|Z\Delta_1\|_{\mathrm{F}}+\frac {1}{\sqrt{N}}\|Z\Delta_2\|_{\mathrm{F}}
\leq J(\|\Delta_1\|_{\mathrm{F}}+\|\Delta_2\|_{\mathrm{F}}) \leq \sqrt{2}J\|\Delta\|_{\mathrm{F}},
\end{align*}
where last inequality follows from H{\"o}lder inequality and $\|\Delta_1\|_{\mathrm{F}}^2+\|\Delta_2\|_{\mathrm{F}}^2= \|\Delta\|_{\mathrm{F}}^2$. Consider that $\|\Delta\|_{\mathrm{F}}\leq \delta\sqrt{\gamma_{\min}(\Sigma)}/(4\sqrt{2})$. Then,
\[
\frac {1}{\sqrt{N}}\|Z\Theta\|_{\mathrm{F}}\leq \frac { J\delta\sqrt{\gamma_{\min}(\Sigma)}}{4}+\sqrt{\gamma_{\max}(\Sigma)}(1+\frac {\delta}{2}).
\]
for any $\Theta\in \mathcal{U}_r $ and hence $J\leq \sqrt{\gamma_{\max}(\Sigma)}(1+\delta/2)/(1-\delta/4)\leq \sqrt{\gamma_{\max}}(1+\delta)$. On the other hand, we have
\begin{align*}
  \frac {1}{\sqrt{N}}\|Z\Theta\|_{\mathrm{F}} &\geq \frac {1}{\sqrt{N}}\|ZP\|_{\mathrm{F}}-\frac {1}{\sqrt{N}}\|Z\Delta\|_{\mathrm{F}} \\
  & \geq  \sqrt{\gamma_{\min}(\Sigma)}(1-\delta/2)  -\sqrt{2}\sqrt{\gamma_{\max}(\Sigma)}(1+\delta)\delta\sqrt{\gamma_{\min}(\Sigma)}/(4\sqrt{2})\\
  &\geq \sqrt{\gamma_{\min}(\Sigma)}(1-\delta).
\end{align*}
This completes the proof. \qed

 In particular, if the measurement matrix $Z$ has entries i.i.d. sampled from a distribution with zero-mean and variance $\sigma^2$, then Theorem \ref{th4} holds with $K_1=K_2=\sigma$, i.e., the RIP holds. Combining with Theorem \ref{th5}, we derive the following corollary directly.

\begin{corollary}\label{col1}
There exist constant $c_1$, $c_2$ such that when $N\geq c_1(2+3s)(m+n)r$, the optimal solution of the SDP program \eqref{noise5} can recover the low-rank matrix $\Theta^\ast$ exactly with probability at least $1-e^{c_2N}$.
\end{corollary}

Corollary \ref{col1} showed that the coefficient matrix of VARX model can be exactly recovered from $\mathcal{O}\left((m+n)r\right)$ samples. Note that an $(n+m)\times n$ matrix of rank $r$ has $r(m+2n-r)$ degrees of freedom. Thereby, this sample size is appropriate.

\subsection{Restricted Operator-norm Curvature and Error Bounds}
 In this subsection, we turn to the noisy setting where we observation $(X,Z)\in \mathbb{R}^{N\times n}\times \mathbb{R}^{N\times (n+m)}$ contacted by $X=Z\Theta+W$. However, it is hard to recover exactly the low-rank matrix in such setting. Here we estimate the low-rank matrix by the program \eqref{m1}. To this end, we first introduce an underlying lemma.
\begin{lemma}\label{lem3} \cite{neg2011}
For the choice of regularization parameter $\lambda_N\geq 2\|\nabla \mathcal{L}(\Theta^\ast)\|_{\mathrm{op}}$, the error matrix $\widehat{\Delta}:=\widehat{\Theta}-\Theta^\ast$ to the program \eqref{m1} satisfies the cone-constrain
\begin{equation}\label{m8}
\|\widehat{\Delta}_{\mathcal{\bar{M}}^{\perp}}\|_{\mathrm{nuc}}\leq 3\|\widehat{\Delta}_{\mathcal{\bar M}}\|_{\mathrm{nuc}}, 
\end{equation}
where $\widehat{\Delta}_{\mathcal{\bar M}}$ has rank at most 2r. Moreover, we have
\begin{equation}\label{m9}
  \|\widehat{\Delta}\|_{\mathrm{nuc}}\leq 4\sqrt{2r}\|\widehat{\Delta}\|_{\mathrm{F}}.
\end{equation}
\end{lemma}
Next, we introduce the notion of restricted curvature that involves a lower bound on the gradient of the cost function.
\begin{definition}\label{def1}
(Operator-norm curvature condition). The cost function $\mathcal{L}_N(\Theta)$ satisfies a \emph{operator-norm curvature condition} with curvature $\mathcal{K}$, tolerance $\tau_N$ if
\[
\|\nabla \mathcal{L}_N(\Theta^\ast+\Delta)-\nabla \mathcal{L}_N(\Theta^\ast)\|_{\mathrm{op}}\geq \mathcal{K}\|\Delta\|_{\mathrm{op}}-\tau_N\|\Delta\|_{\mathrm{nuc}}
\]
for all $\Delta\in \mathbb{R}^{(n+m)\times n}$.
\end{definition}
The definition is easier to understand and apply than RSC. Under the operator-norm curvature condition, we can prove the following theorem directly.
\begin{theorem}\label{th1}
Suppose the cost function $\mathcal{L}_N(\Theta)$ satisfies the operator-norm curvature condition with parameters $(\mathcal{K},\tau_N)$, and the true matrix $\Theta^\ast$ with rank $r<\frac {\mathcal{K}}{64\tau_N}$. Then when the regularization parameter $\lambda_N\geq 2\|\nabla\mathcal{L}(\Theta^\ast)\|_{\emph{op}}$ the estimator $\widehat{\Theta}$ to the program \eqref{m1} satisfies the bound
\[
\|\widehat{\Theta}-\Theta^\ast\|_{\emph{op}}\leq 3\frac {\lambda_N}{\mathcal{K}}.
\]
\end{theorem}
\noindent {\bf Proof.}\  For the SDP program \eqref{m1}, by the KKT condition \cite{boy2004}, there exists $G\in \partial \|\widehat{\Theta}\|_{\mathrm{nuc}}$ with $\|G\|_{\mathrm{op}}\leq1$ such that
\[
\nabla \mathcal{L}(\widehat{\Theta})+\lambda_N G=0,
\]
where $\partial\|\widehat{\Theta}\|_{\mathrm{nuc}}$ stands for the subdifferential \cite{wat1992} of the nuclear norm at $\widehat{\Theta}$. Define $\Delta:=\widehat{\Theta}-\Theta^\ast$.
By the triangle inequality and the inequality \eqref{m8}, when $\lambda_N\geq 2\|\nabla\mathcal{L}(\Theta^\ast)\|_{\mathrm{op}}$, one has
\[
\|\nabla \mathcal{L}(\Theta^\ast+\widehat{\Delta})-\nabla \mathcal{L}(\Theta^\ast)\|_{\mathrm{op}}\leq \|\nabla \mathcal{L}(\Theta^\ast)\|_{\mathrm{op}}+\lambda_N.
\]
Applying the operator-norm curvature condition yields
\begin{equation}\label{m10}
 \mathcal{K}\|\widehat{\Delta}\|_{\mathrm{\mathrm{op}}}\leq \tau_N\|\widehat{\Delta}\|_{\mathrm{nuc}}+\frac {3}{2}\lambda_N.
\end{equation}
Moreover, by H\"older's inequality and \eqref{m9}, one has
\[
\|\widehat{\Delta}\|_{\mathrm{F}}^2 \leq \|\widehat{\Delta}\|_{\mathrm{nuc}}\|\widehat{\Delta}\|_{\mathrm{op}}\leq 4\sqrt{2r}\|\widehat{\Delta}\|_{\mathrm{F}} \|\widehat{\Delta}\|_{\mathrm{op}},
\]
and hence $\|\widehat{\Delta}\|_{\mathrm{F}}\leq 4\sqrt{2r}\|\widehat{\Delta}\|_{\mathrm{op}}$. Bringing it into \eqref{m9} yields
\begin{equation}\label{m11}
   \|\widehat{\Delta}\|_{\mathrm{nuc}}\leq 32r \|\widehat{\Delta}\|_{\mathrm{op}}.
\end{equation}
Finally, substituting \eqref{m11} into \eqref{m10} yields the desired claim. \qed
\begin{remark}
It is noting that the tolerance $\tau_N$ can be zero. In this case, the upper bound constraint for rank can be neglected. Moreover, the operator-norm bound in Theorem \ref{th1} does not depend on $r$. However, the Frobenius norm bound has the scaling parameter $\sqrt{r}$ \cite{neg2011,neg2009}. In some sense, the operator-norm bound is stronger than the Frobenius norm bound. This theorem will then be used to obtain a non-asymptotic result.

\end{remark}
Next, we turn to the identification of dynamic system \eqref{p1}.
  Observe that
\begin{align}\notag
  x^{(i)}(T^0-1)&= A^{T^0-2}Bu^{(i)}(0)+\cdots+Bu^{(i)}(T^0-2) \\ \label{m2}
   &+A^{T^0-2}w^{(i)}(0)+\cdots+w^{(i)}(T^0-2).
 \end{align}
Assume that the noise vector $w(t)$ is sub-Gaussian with parameter $\sigma_w$. Using sub-Gaussian excitation with parameter $\sigma_u$, which is independent of $w(t)$, we then have the following results.

 \begin{proposition}\label{pro1}
 The random vector $z^{(i)}(T^0-1)\in \mathbb{R}^{n+m}$ (the $i^{th}$ row of $Z$) is sub-Gaussian with parameter $\sigma_z$ for any $i=1,2,\cdots,N$, where $\sigma_z=\sqrt{\sum \limits_{k=0}^{T^0-2}\left(\|A^{T^0-2-k}B\|_{\emph{op}}\sigma_u^2+\|A^{T^0-2-k}\|_{\emph{op}}\sigma_w^2\right)+\sigma_u^2}$.
 \end{proposition}
\noindent {\bf Proof.}\ For any fixed $v\in \mathbb{R}^{n+m}$ and any $\lambda \in \mathbb{R}$, by Eq. \eqref{m2} and the independence of $u^{(i)}(t)$ and $w^{(i)}(t)$, $t=1,\cdots,T^0-1$, we have
\begin{align*}
  \mathbb{E} e^{\lambda \langle v,z^{(i)}(T^0-1)\rangle}=&\mathbb{E}e^{\lambda \langle v_1,x^{(i)}(T^0-1) \rangle}\mathbb{E} e^{\lambda \langle v_2,u^{(i)}(T^0-1) \rangle} \\
  =& \mathop {\prod} \limits_{k=1}^{T^0-2}\left(\mathbb{E} e^{\lambda \langle v_1,A^{T^0-2-k}Bu^{(i)}(k) \rangle} \mathbb{E} e^{\lambda \langle v_1,A^{T^0-2-k}w^{(i)}(k) \rangle}\right)\mathbb{E}e^{\lambda \langle v_2,u^{(i)}(T^0-1) \rangle}\\
  =& \mathop {\prod} \limits_{k=1}^{T^0-2} \left(\mathbb{E} e^{\lambda \langle (A^{T^0-2-k}B)^Tv_1,u^{(i)}(k) \rangle} \mathbb{E} e^{\lambda \langle (A^{T^0-2-k})^Tv_1,w^{(i)}(k) \rangle}\right)\mathbb{E}e^{\lambda \langle v_2,u^{(i)}(T^0-1) \rangle}\\
   \leq& \mathop {\prod} \limits_{k=1}^{T^0-2} \left(e^{\frac {\lambda^2 \|(A^{T^0-2-k}B)^Tv_1\|_2^2\sigma_u^2+\lambda^2 \|(A^{T^0-2-k})^Tv_1\|_2^2\sigma_w^2}{2}}   \right)e^{\frac {\lambda^2\|v_2\|_2^2\sigma_u^2}{2}} \\
  \leq & e^{\frac{\lambda^2\sigma_z^2}{2}},
\end{align*}
where $v=\left[
           \begin{array}{c}
             v_1 \\
             v_2 \\
           \end{array}
         \right]
$, $v_1\in \mathbb{R}^n$, $v_2\in \mathbb{R}^m$, which implies that $z^{(i)}(T^0-1)$ is sub-Gaussian with parameter $\sigma_z$.    \qed

Proposition \ref{pro1} showed that the rows of the measurement matrices are i.i.d. sub-Gaussian.        Denote by $\Sigma$ the covariance matrix of random vector $z^{(i)}(T^0-1)$. Assume that $\Sigma$ is invertible.
Applying the triangle inequality and \eqref{m3} yields
\begin{eqnarray*}
  \|\widehat{\Sigma}\|_{\mathrm{op}} &\leq & \|\widehat{\Sigma}-\Sigma\|_{\mathrm{op}}+\|\Sigma\|_{\mathrm{op}}
 \\
  & \leq &  16\sqrt{6}\beta^2\left(\sqrt {\frac {n+m}{N}}+\frac {n+m}{N}\right)+\delta\beta^2 +\gamma_{\max}\\
  & \leq & ( 32\sqrt{6}+\delta)\beta^2 +\gamma_{\max}(\Sigma)
\end{eqnarray*}
with probability at least $1-e^{-N \min\{\frac{\delta}{16\sqrt{2}},\frac{\delta^2}{512}\}}$, where $\gamma_{\max}(\Sigma)$ represents the maximal singular value of the covariance matrix $\Sigma$. In particular, let $\delta=1$, we have
\begin{equation}\label{m4}
  \mathbb{P}\left[ \|\widehat{\Sigma}\|_{\mathrm{op}}\geq ( 32\sqrt{6}+1)\beta^2 +\gamma_{\max}(\Sigma) \right]\leq e^{-\frac{N}{512}}.
\end{equation}

\begin{lemma}\label{lem2}
 There exist non-negative constants $c_1$, $c_2$, such that
\[
\mathbb{P}\left[ \left\|\frac {Z^T W}{N}\right\|_{\mathrm{op}}\geq 2\alpha\sqrt \frac {(n+m)}{N} \right]\leq c_1e^{-c_2(n+m)},
\]
where $\alpha^2:=2\sigma_w^2[( 32\sqrt{6}+1)\beta^2 +\gamma_{\max}(\Sigma)]$.
\end{lemma}
\noindent {\bf Proof.}\  See Appendix \ref{app3}. \qed

Applying the Theorem \ref{th1}, we obtain a operator norm error bound, which holds with high probability.

\begin{corollary}\label{th2} Consider the multivariate system \eqref{p1} where true parameter matrix $\Theta^\ast$ has rank $r<n$. There are constants $c_1$, $c_5$, $c_6$ such that when $N\geq \frac {c_1\beta^4}{\gamma_{\min}^2(\Sigma)}(n+m)$ the solution to the program \eqref{m1} with $\lambda_N=4\alpha\sqrt \frac {n+m}{N}$ satisfies the bound
\begin{equation}\label{m6}
\|\widehat{\Theta}-\Theta^\ast\|_{\mathrm{op}}\leq   \frac {12\alpha}{\gamma_{\min}(\Sigma)} \sqrt{\frac {n+m}{N}}
\end{equation}
with probability at least $1-c_5e^{-c_6(n+m)}$, where $\alpha^2:=2\sigma_w^2[( 32\sqrt{6}+1)\beta^2 +\gamma_{\max}(\Sigma)]$.
\end{corollary}
\noindent {\bf Proof.} We first prove that the curvature condition holds for the parameters $\mathcal{K}=\frac {\gamma_{\min}(\Sigma)}{2}$, $\tau_N=0$ with high probability.
  Since $\nabla\mathcal{L}_N(\Theta)=\frac {1}{N}Z^T(X-Z\Theta)$, then for any $\Delta \in\mathbb{R}^{(n+m)\times n} $, we have $\nabla\mathcal{L}_N(\Theta^\ast+\Delta)-\nabla\mathcal{L}_N(\Theta^\ast)=-\widehat{\Sigma}\Delta$, where $\widehat{\Sigma}=\frac {1}{N}Z^T Z$ stands for the sample covariance matrix. On the other hand, one has
\begin{eqnarray*}
\| \widehat{\Sigma}\Delta\|_{\mathrm{op}} &\geq& \| \Sigma \Delta\|_{\mathrm{op}}-\| (\widehat{\Sigma}-\Sigma)\Delta\|_{\mathrm{op}}  \\
   &\geq& \gamma_{\min}(\Sigma)\|\Delta\|_{\mathrm{op}}-\|\widehat{\Sigma}-\Sigma\|_{\mathrm{op}}\|\Delta\|_{\mathrm{op}}.
\end{eqnarray*}
By Lemma \ref{lem1}, there exists a constant $\beta$ such that
\begin{equation*}
                              \| \widehat{\Sigma}\Delta\|_{\mathrm{op}}
                               \geq  \gamma_{\min}(\Sigma)\|\Delta\|_{\mathrm{op}}-\left(16\sqrt{6}\left(\sqrt {\frac {n+m}{N}}+\frac {n+m}{N}\right)+\delta\right)\beta^2  \|\Delta\|_{\mathrm{op}}
                            \end{equation*}
with probability greater than $1-e^{-N \min\big{\{}\frac{\delta}{16\sqrt{2}},\frac{\delta^2}{512}\big{\}}}$.  Taking $\delta=\frac {\gamma_{\min}(\Sigma)}{4\beta^2}$, there exist positive constants $c_1,c_2$ such that when $N\geq \frac {c_1\beta^4}{\gamma_{\min}^2(\Sigma)}(n+m)$,
\[
\| \widehat{\Sigma}\Delta\|_{\mathrm{op}}\geq \frac {1}{2}\gamma_{\min}(\Sigma)\|\Delta\|_{\mathrm{op}}
\]
with probability greater than $1-e^{-c_2 N }$, where we use the fact that
\[
\sqrt{\frac {n+m}{N}}+\frac {n+m}{N}\leq 2 \sqrt{\frac {n+m}{N}},\ \ \mathrm{when} \ N\geq n+m.
\]
Thus, the operator-norm curvature condition holds with curvature $\mathcal{K}=\gamma_{\min}(\Sigma)$ and tolerance $\tau_N=0$ . On the other hand, we prove that the inequality $\|\nabla\mathcal{L}(\Theta^\ast)\|_{\mathrm{op}}\leq \frac {\lambda_N}{2}$ holds with high probability for the selected regularization parameter $\lambda_N$. Consider that $\nabla\mathcal{L}(\Theta^\ast)=\frac {Z^T W}{N}$. Since the rows of matrix $W$ are i.i.d. sub-Gaussian vectors, then by Lemma \ref{lem2}, there exist constants $c_3,c_4$ such that
\[
\mathbb{P}\left[ \left\|\frac {Z^T W}{N}\right\|_{\mathrm{op}}\geq 2\alpha\sqrt \frac {n+m}{N} \right]\leq c_3e^{-c_4(n+m)},
\]
or equivalently,
\[
\mathbb{P}\left[\|\nabla\mathcal{L}(\Theta^\ast)\|_{\mathrm{op}}\geq \frac {\lambda_N}{2}\right]\leq c_3e^{-c_4(n+m)}.
\]
Applying Theorem \ref{th1} yields
\[
\|\widehat{\Theta}-\Theta^\ast\|_{\mathrm{op}}\leq  \frac {24\alpha}{\gamma_{\min}(\Sigma)} \sqrt{\frac {n+m}{N}}
\]
with probability at least $1-e^{-c_2 N }- c_3e^{-c_4(n+m)}$, which implies \eqref{m6}.

\begin{remark}\label{remark2}
By \eqref{m9}, we have
\begin{align*}
  \|\widehat{\Theta}-\Theta^\ast\|_{\mathrm{F}}^2 & \leq \|\widehat{\Theta}-\Theta^\ast\|_{\mathrm{nuc}}\|\widehat{\Theta}-\Theta^\ast\|_{\mathrm{op}} \\
   & \leq 4\sqrt{2r}\|\widehat{\Theta}-\Theta^\ast\|_{\mathrm{F}}\|\widehat{\Theta}-\Theta^\ast\|_{\mathrm{op}}.
\end{align*}
Applying Corollary \ref{th2} yields
\begin{equation}\label{m14}
  \|\widehat{\Theta}-\Theta^\ast\|_{\mathrm{F}}\leq   \frac {96\sqrt{2r}\alpha}{\gamma_{\min}(\Sigma)} \sqrt{\frac {n+m}{N}}
\end{equation}
with probability at least $1-c_5e^{-c_6(n+m)}$.
\end{remark}
In \cite{neg2011}, the author discusses the RSC condition, under which the nonasymptotic error bound is obtained with high probability. Here, by Lemmas \ref{lem1} and \ref{lem2}, we can derive a Frobenius norm bound like \eqref{m14}. Obviously, the operator norm bound \eqref{m6} is stronger than the Frobenius norm bound in some sense. Moreover, this samples size has smaller scaling parameter $m+n$ than the scaling parameter $r(m+n)$ obtained in many studies.

The above results are stated for matrices that are exactly low rank. Moreover,  there are weak low-rank matrices that can be closely approximated by low rank matrix. For a parameter $q\in [0,1]$ and radius $R_q>0$, consider the set
\[
\mathbb{B}(R_q)=\Big{\{}\Theta \in \mathbb{R}^{(n+m)\times n} \big{|}\sum _{j=1}^{n}|\sigma_j(\Theta)|^{q}\leq R_q \Big{\}}.
\]
In the special case $q=0$, any $ \Theta^\ast\in \mathbb{B}(R_0)$ is a matrix whose rank not exceed $R_q$. Then we have the following corollary.
\begin{theorem}\label{th7}
Suppose the true matrix $\Theta^\ast\in \mathbb{B}(R_q)$, the regularization parameter $\lambda_N\geq 2\|\nabla\mathcal{L}(\Theta^\ast)\|_{\emph{op}}$, and the cost function $\mathcal{L}_N(\Theta)$ satisfies the operator-norm curvature condition with parameters $(\mathcal{K},\tau_N)$. Then when $R_q\leq \frac {\mathcal{K}}{128\tau_N}$, the estimator $\widehat{\Theta}$ to the program \eqref{m1} satisfies the bound
\[
\|\widehat{\Theta}-\Theta^\ast\|_{\emph{op}}\leq \max  \Big{ \{ }  \frac {32\tau_N R_q}{\mathcal{K}}, \frac {6}{\mathcal{K}}\lambda_N \Big{ \}}.
\]
\end{theorem}
\noindent {\bf Proof.}\
The matrix $\Theta^\ast$ can be decomposed as $\Theta^\ast=UDV^T$, where $U\in \mathbb{R}^{(n+m)\times n}$ and $V\in \mathbb{R}^{n\times n}$ are orthogonal matrices, and $D$ is a diagonal with its entries corresponding to the singular values in non-increasing order $\sigma_1(\Theta^\ast)\geq \sigma_2(\Theta^\ast)\cdots \geq \sigma_n(\Theta^\ast)\geq 0$. Define
\[
S:=\big{\{}j\in \{1,2\cdots, n\}| \sigma_j(\Theta^\ast)> \tau  \big{\}}
\]
with parameter $\tau>0$. Define
\[
\Theta^\prime:=U\left[
                 \begin{array}{cc}
                   0 & 0 \\
                   0 & D_1 \\
                 \end{array}
               \right]V^T,
\]
where $D_1:=\mathrm{diag}\left(\sigma_{(|S|+1)}(\Theta^\ast),\cdots,\sigma_n(\Theta^\ast)\right)$. Then, $\mathrm{rank}(\Theta^\prime)\leq n-|S|$ and
\begin{equation}\label{ne1}
  \|\Theta^\prime\|_{\mathrm{nuc}}=\sum \limits_{j=|S|+1}^{n}\sigma_j(\Theta^\ast)\leq \tau \sum \limits_{j=|S|+1}^{n}\left(\sigma_j(\Theta^\ast)/\tau \right)^q=\tau^{1-q}R_q.
\end{equation}
Moreover, by definition of  $\mathbb{B}(R_q)$, one has
\begin{equation}\label{ne2}
|S|\leq  \sum \limits_{j=1}^{n}\left(\sigma_j(\Theta^\ast)/\tau \right)^q=\tau^{-q}R_q.
\end{equation}
Define $\widehat{\Delta}:=\widehat{\Theta}-\Theta^\ast$. By the proof of Theorem \ref{th1}, one has
\begin{equation}\label{ne3}
 \mathcal{K}\|\widehat{\Delta}\|_{\mathrm{\mathrm{op}}}\leq \tau_N\|\widehat{\Delta}\|_{\mathrm{nuc}}+\frac {3}{2}\lambda_N.
\end{equation}
By Lemma 1 in \cite{neg2011}, there exists a matrix decomposition, $\widehat{\Delta}=\Delta_1+\Delta_2$, where $\mathrm{rank}(\Delta_1)\leq |S|$ and
\begin{equation}\label{ne4}
  \|\Delta_2\|_{\mathrm{nuc}}\leq 3\|\Delta_1\|_{\mathrm{nuc}}+4\|\Theta^\prime\|_{\mathrm{nuc}}.
\end{equation}
Consider that $\|\Delta_1\|_{\mathrm{nuc}}\leq \sqrt{2|S|}\|\Delta_1\|_{\mathrm{F}}$. Combining with the inequality \eqref{ne4} yields
  \begin{align*}
  \|\widehat{\Delta}\|_{\mathrm{F}}^2  &  \leq \|\widehat{\Delta}\|_{\mathrm{nuc}}\|\widehat{\Delta}\|_{\mathrm{op}}\leq \left(4\|\Delta_1\|_{\mathrm{nuc}}+4\|\Theta^\prime\|_{\mathrm{nuc}}\right)\|\widehat{\Delta}\|_{\mathrm{op}}\\
  &\leq  4\left(\sqrt{2|S|}\|\widehat{\Delta}\|_{\mathrm{F}}+\|\Theta^\prime\|_{\mathrm{nuc}}\right)\|\widehat{\Delta}\|_{\mathrm{op}}.                                                       \end{align*}
It leads to
\begin{equation}\label{ne6}
\|\widehat{\Delta}\|_{\mathrm{F}}\leq \max  \Big{\{} 8\sqrt{2|S|}\|\widehat{\Delta}\|_{\mathrm{op}},2\sqrt{2\|\Theta^\prime\|_{\mathrm{nuc}} \|\widehat{\Delta}\|_{\mathrm{op}} } \Big{\}}.
\end{equation}
Then, by the inequalities \eqref{ne2}, \eqref{ne4} and \eqref{ne6}, one has
\begin{align*}
    \|\widehat{\Delta}\|_{\mathrm{nuc}} & \leq  \|\Delta_1\|_{\mathrm{nuc}}+ \|\Delta_2\|_{\mathrm{nuc}}\leq 4\sqrt{2|S|}\|\widehat{\Delta}\|_{\mathrm{F}} +4\|\Theta^\prime\|_{\mathrm{nuc}}\\
                                        &  \leq 64|S|\|\widehat{\Delta}\|_{\mathrm{op}}+16\sqrt{|S|\|\Theta^\prime\|_{\mathrm{nuc}} \|\widehat{\Delta}\|_{\mathrm{op}} }+4\|\Theta^\prime\|_{\mathrm{nuc}}\\
                                        & \leq 64\tau^{-q}R_q\|\widehat{\Delta}\|_{\mathrm{op}}+16\sqrt{\tau^{-q}R_q\|\Theta^\prime\|_{\mathrm{nuc}} \|\widehat{\Delta}\|_{\mathrm{op}} }+4\|\Theta^\prime\|_{\mathrm{nuc}}
\end{align*}\
Substituting this inequality into the inequality \eqref{ne3} yields
\[
 \left(1-\frac {64\tau_N\tau^{-q}R_q}{\mathcal{K}}\right)\|\widehat{\Delta}\|_{\mathrm{\mathrm{op}}}\leq \frac {16\tau_N  \sqrt{\tau^{-q}R_q\|\Theta^\prime\|_{\mathrm{nuc}} \|\widehat{\Delta}\|_{\mathrm{op}} } }{\mathcal{K}}+\frac {3}{2\mathcal{K}}\lambda_N.
\]
Then, when $R_q\leq \frac {\mathcal{K}\tau^{q}}{128\tau_N}$, one has
\[
 \|\widehat{\Delta}\|_{\mathrm{\mathrm{op}}}\leq 2\sqrt{\frac {2\tau_N\|\Theta^\prime\|_{\mathrm{nuc}} \|\widehat{\Delta}\|_{\mathrm{op}} }{\mathcal{K}}}+\frac {3}{\mathcal{K}}\lambda_N.
\]
Hence,
\[
\|\widehat{\Delta}\|_{\mathrm{\mathrm{op}}}\leq  \max  \Big{ \{ }  \frac {32\tau_N\tau^{1-q}R_q}{\mathcal{K}}, \frac {6}{\mathcal{K}}\lambda_N \Big{ \}}.
\]
Letting $\tau=1$, we can obtain the desired results. \qed

 Based on Corollary \ref{th2} and Theorem \ref{th7}, we can obtain the following result directly.
\begin{corollary}\label{col5}
Consider the multivariate system \eqref{p1} where true parameter matrix $\Theta^\ast \in\mathbb{B}(R_q)$. There are constants $c_1$, $c_5$, $c_6$ such that when $N\geq \frac {c_1\beta^4}{\gamma_{\min}^2(\Sigma)}(n+m)$ the solution to the program \eqref{m1} with $\lambda_N=4\alpha\sqrt \frac {n+m}{N}$ satisfies the bound
\begin{equation}\label{m6}
\|\widehat{\Theta}-\Theta^\ast\|_{\mathrm{op}}\leq   \frac {24\alpha}{\gamma_{\min}(\Sigma)} \sqrt{\frac {n+m}{N}}
\end{equation}
with probability at least $1-c_5e^{-c_6(n+m)}$, where $\alpha^2:=2\sigma_w^2[( 32\sqrt{6}+1)\beta^2 +\gamma_{\max}(\Sigma)]$.
\end{corollary}

Corollary \ref{col5} showed that the selected regularization parameter $\lambda_N$ is decreasing as the sample size $N$ increases.  In particular, if $q=0$, the above results for weak low rank matrices can be boiled down to the results for exactly low rank matrices. Clearly, this operator norm error bound is a stronger bound than the common Frobenius norm error bound.

\section{Discussions}\label{conclusion}

The paper is concerned with the identification of dynamic systems with low-rank constraints. Given a finite number of sample sets, the aim is to recover the coefficient matrix of the dynamic systems. For this purpose, nuclear norm heuristic is considered. This paper focuses on a VARX(1) model. By the sampling method and the more general input designs, we generalized the existing work. Moreover, it can be generalized to VARX(d) model. The form of VARX(d) model is as follows:
\[
x(t+1)=\sum \limits_{k=0}^{d-1}A_kx_{t-k}+Bu(t)+w(t).
\]
This model can be rendered into the multivariate regression model \eqref{int1} using the the standard transformation
\[
\left[
   \begin{array}{c}
     x(t+1) \\
     0\\
      \vdots \\
       0\\
   \end{array}
 \right]
 =\left[
           \begin{array}{cccc}
             A_0 & A_1 & \cdots & A_{d-1} \\
             0 & 0 & \cdots & 0 \\
             \vdots & \vdots &  & \vdots \\
             0 & 0 & \cdots & 0 \\
           \end{array}
         \right] \left[
                   \begin{array}{c}
                     x(t) \\
                     x({t-1}) \\
                      \vdots\\
                       x({t-d+1})\\
                   \end{array}
                 \right]+
                 \left[
                   \begin{array}{c}
                     B \\
                      0\\
                     \vdots \\
                     0 \\
                   \end{array}
                 \right]u(t) +
                 \left[
                   \begin{array}{c}
                     w(t) \\
                     0 \\
                    0 \\
                     0 \\
                   \end{array}
                 \right]
\]
Note that, the conclusions in this paper are also applicable to this problem. Further, our research will focus on the dynamics systems with low-rank and sparse structure simultaneously.

%

\appendix

\section{Proof of Lemma~\ref{lem2}}\label{app3}

Let $\mathcal{A}:=\{u^1,u^2,\cdots,u^J\}$ and $\mathcal{B}:=\{v^1,v^2,\cdots,v^L\}$ be 1/4 coverings of the spheres $\mathbb{S}^{n+m-1}$ and  $\mathbb{S}^{n-1}$, respectively. By the spherical covering theorem \cite{ver2010}, there exists a covering of the spheres $\mathbb{S}^{n+m-1}$ and  $\mathbb{S}^{n-1}$ with $J\leq 9^{n+m}$, $L \leq 9^n$ elements respectively. For any $v\in \mathbb{S}^{n-1}$, there exist a vector $v^{l}\in \mathcal{B}$ such that $v=v^l+v^{\prime}$ with $\|v^{\prime}\|_2\leq 1/4$. Define $Q:=\frac {Z^TW}{N}$. Then
\[
\|Q\|_{\mathrm{op}}=\sup \limits_{v\in \mathbb{S}^{n-1}}\|Qv\|_{2}\leq \sup \limits_{l=1,2,\cdots,L}\|Qv^l\|_{2}+\frac {1}{4}\|Q\|_{\mathrm{op}}.
\]
Next, by similar argument, we have $\|Qv^l\|_{2}\leq \sup \limits_{j=1,2,\cdots,J}|\langle u^j,Qv^l\rangle|+\frac {1}{4}\|Q\|_{\mathrm{op}}$ for any fixed $l$. Thus,
\begin{align}\notag
\|Q\|_{\mathrm{op}}&\leq 2\sup \limits_{j=1,2,\cdots,J}\sup \limits_{l=1,2,\cdots,L}|\langle u^j,Qv^l\rangle|\\\label{m5}
&=2\sup \limits_{j=1,2,\cdots,J}\sup \limits_{l=1,2,\cdots,L}\frac {1}{N}\sum \limits_{k=1}^{N}\langle u^j,Z_{k,:}^T\rangle\langle W_{k,:}^T,v^l\rangle. \tag{A.1}
\end{align}
Since the rows of noise matrix $W$ are independent of each other and obey the sub-Gaussian distribution with parameter $\sigma_w$, $\langle W_{k,:}^T,v^l\rangle\ (k=1,\cdots,N)$ are the sub-Gaussian variables with parameter $\sigma_w$. Besides, the random matrices $Z$ and $W$ are independent. Therefore, conditioned on $Z$, the variable $U:=\frac {1}{N}\sum \limits_{k=1}^{N}\langle u^j,Z_{k,:}^T\rangle\langle W_{k,:}^T,v^l\rangle$ is sub-Gaussian with parameter $\tilde{\sigma}:=\frac {\sigma_w}{N}\|Zu^j\|_2$. Define the event $\mathcal{F}=\{\tilde{\sigma}^2\leq \frac {\sigma_w^2}{N}[( 32\sqrt{6}+1)\beta^2 +\gamma_{\max}(\Sigma)]\}$. Consider that
\[
\tilde{\sigma}^2=\frac {\sigma_w^2}{N}(u^j)^T\widehat{\Sigma} u^j\leq \frac {\sigma_w^2}{N}\|\widehat{\Sigma}\|_{\mathrm{op}}.
\]
Then by \eqref{m4}, we have $\|\widehat{\Sigma}\|_{\mathrm{op}}\leq ( 32\sqrt{6}+1)\beta^2 +\gamma_{\max}(\Sigma)$ with probability at least $1-e^{-\frac{N}{512}}$, and hence $\mathbb{P}\left[\mathcal{F}^c\right] \leq e^{-\frac{N}{512}}$. Moreover, by the total probability rule and the sub-Gaussian tail bounds, we have
\begin{eqnarray*}
  \mathbb{P}\left[ |U|\geq t\right] & \leq & \mathbb{P}\left[ |U|\geq t|\mathcal{F}\right]+\mathbb{P}\left[ \mathcal{F}^c\right]\\
   &\leq &  2e^{-N\frac {t^2}{\alpha^2}}+e^{-\frac{N}{512}},
\end{eqnarray*}
where $\alpha^2:=2\sigma_w^2[( 32\sqrt{6}+1)\beta^2 +\gamma_{\max}(\Sigma)]$. Therefore, by \eqref{m5}, we have
\[
\mathbb{P}\left[ \|Q\|_{\mathrm{op}}\geq 2t \right]\leq 9^{m+2n} ( 2e^{-N\frac {t^2}{\alpha^2}}+e^{-\frac{N}{512}}).
\]
Setting $t^2= \frac {(n+m)\alpha^2}{512N}$, we establish the desired claim. \qed

\end{document}